\documentclass[11pt,reqno]{amsproc} 
\usepackage{amsmath}
\usepackage{amssymb} 
\usepackage[all]{xy} 
\usepackage{psfrag}
\usepackage{graphicx}

\theoremstyle{plain}
\newtheorem{theorem}{Theorem}

\newtheorem{lemma}[theorem]{Lemma}

\theoremstyle{definition}

\newtheorem{example}[theorem]{Example}
\newtheorem{remark}[theorem]{Remark}
\newtheorem{notation}[theorem]{Notation}

\newcommand{\norm}[2]{\left\| {#1}\right\| _{#2}}
\newcommand{\abs}[1]{\left |{#1}\right |}

\renewcommand{\l}{{\mathcal L}} 
 
\newcommand{\set}[1]{\left\{{#1}\right \}}

\newcommand{\all}[2]{\{ {#1}\,|\,{#2} \}}

\newcommand{\I}{\mathcal I}
\newcommand{\C}{\mathbb C}

\newcommand{\N}{\mathbb N}

\newcommand{\un}{\underline{n}}
\newcommand{\at}{{\tilde{\alpha}}}
\newcommand{\bt}{{\tilde{\beta}}}

\title[Transfer operator eigenvalues]{Explicit a priori bounds on
transfer operator eigenvalues}

\author{Oscar F.~Bandtlow and Oliver Jenkinson}

\address{Oscar F.~Bandtlow; 
School of Mathematical Sciences, Queen
  Mary, University of London, Mile End Road, London, E1 4NS, UK.
\newline
{\tt ob@maths.qmul.ac.uk} \newline {\tt
    www.maths.qmul.ac.uk/$\sim$ob}}  

\address{Oliver Jenkinson; 
School of Mathematical Sciences, Queen
  Mary, University of London, Mile End Road, London, E1 4NS, UK.
  \newline {\tt omj@maths.qmul.ac.uk} \newline {\tt
    www.maths.qmul.ac.uk/$\sim$omj}}

\begin{document}

\date{\today}

\begin{abstract}
We provide explicit bounds on the eigenvalues of 
transfer operators defined in terms of holomorphic data.
\end{abstract}

\maketitle

Linear operators of the form
$\l f= \sum_{i\in\I} w_i\cdot f\circ T_i$,
so-called \emph{transfer operators} 
(see e.g.~\cite{baladibook, ruellecmp68,ruelleinventiones}),
arise in a number
of problems in 
dynamical systems.
If the $T_i$ are inverse branches of an expanding map $T$, and the
weight functions $w_i$ are positive, the spectrum of $\l$ has well-known
interpretations in terms of the exponential mixing rate  
of an invariant Gibbs measure (see \cite{baladibook}).
Applications also arise when the $w_i$ are real-valued 
(e.g.~\cite{ccr,jms,pollicott})
or complex-valued (e.g.~\cite{dolgopyat,pollicottsharp}).

In this article
we suppose that
$T_i$ and $w_i$ 
are analytic functions of $d$ variables,
for each $i$ in 
some countable\footnote{Subsequent results are new even when 
$\I$ is finite, but it is convenient
to also allow countably infinite $\I$.}
 index set $\I$. 
Under suitable hypotheses on $T_i$ and $w_i$ the 
transfer operator $\l$ defines
a compact operator on Hardy space
$H^2(B)$, and
we can give 
completely explicit bounds on its eigenvalue 
sequence\footnote{Precisely,
$\set{\lambda_n(\l)}_{n=1}^\infty$ denotes the 
sequence of all
eigenvalues of $\l$ counting algebraic
multiplicities and ordered by decreasing modulus,
with the
usual convention (see e.g.~\cite[3.2.20]{pietsch}) that 
distinct eigenvalues with the same modulus can be written in any order.}
$\{\lambda_n(\l)\}_{n=1}^\infty$:

\begin{theorem}
\label{mainresult}
Suppose there is
a complex Euclidean ball $B\subset \C^d$ such that
each $w_i:B\to\C$ is holomorphic with
$\sum_{i\in\I}\sup_{z\in B}|w_i(z)|<\infty$,
and each $T_i:B\to B$ is holomorphic with
$\cup_{i\in\I} T_i(B)$ contained in the
ball concentric with $B$ whose radius is $r<1$
times that of $B$.

Then $\l:H^2(B)\to H^2(B)$ is compact and 
\begin{equation}
\label{generalevalue}
|\lambda_n(\l)| <
\frac{W \sqrt{d}}{r^d (1-r^2)^{d/2}}\ n^{(d-1)/(2d)}\ 
r^{\frac{d}{d+1}(d!)^{1/d}n^{1/d}}\quad\text{for all }n\ge1 \,,
\end{equation}
where $W:=\sup_{z\in B}\sum_{i\in\I}|w_i(z)|$. 

If $d=1$ then
\begin{equation}
\label{d=1evalue}
 |\lambda_n(\l)|\leq \frac{W}{\sqrt{1-r^2}}\  r^{(n-1)/2}
\quad\text{for all }n\ge1 \,.
\end{equation}

\end{theorem}

\begin{remark}
\item[\, (i)]
An estimate of the form
$
|\lambda_n(\l)|\le C\theta^{n^{1/d}}
$
for some (undefined) constants $C>0$, $\theta\in(0,1)$
is asserted, either implicitly or explicitly, in the work of
several authors
(e.g.~\cite{faureroy,fried,glz});
the novelty here is that careful derivation of this bound
renders explicit the constants $C$, $\theta$.
\item[\, (ii)]
Using different techniques, the bound
$
|\lambda_n(\l)|\le C\theta^{n^{1/d}}
$
can also be established in the case where $B$ is an arbitrary
open subset of $\C^d$ (see \cite{explicit}), though here our expressions
for $C$, $\theta$ are more complicated.
\end{remark}

\begin{example}
If 
$\l f(z) = \sum_{n=1}^\infty \left( \frac{1}{n+z}\right)^{2}
f\left(\frac{1}{n+z}\right)$
(the Perron-Frobenius operator 
for the Gauss
map
$x\mapsto 1/x \pmod 1$, cf.~\cite{mayercmp}),
$B\subset\C$ may be chosen as 
the open disc of radius $3/2$ centred at the point $1$.
In this case
$W=
\sup_{z\in B}  
\sum_{n=1}^\infty 
|n+z|^{-2}
=\sum_{n=1}^\infty (n-1/2)^{-2}=\pi^2/2$
and $r=2/3$,
so (\ref{d=1evalue})
yields
\begin{equation*}
 |\lambda_n(\l)|\leq 
\frac{3\pi^2}{2\sqrt{5}}   (2/3)^{(n-1)/2}
\quad\text{for all }n\ge1 \,.
\end{equation*}
\end{example}

\begin{notation}  
For an open ball $D\subset\C^d$, let
$H^\infty(D)$ denote the Banach space 
consisting of all bounded
holomorphic $\C$-valued functions on $D$, 
with norm
$ \norm{f}{H^\infty(D)}:=\sup_{z\in D}|f(z)|$.

\textit{Hardy space}  $H^2(D)$ 
 (see \cite[Ch.~8.3]{krantz})
is the
$L^2(\partial D,\sigma)$-closure of 
the set of those $f\in H^\infty(D)$ which extend continuously to the 
boundary $\partial D$, 
where $\sigma$ denotes
$(2d-1)$-dimensional Lebesgue
measure on $\partial D$,
normalised
so that $\sigma(\partial D)=1$.
In particular, $H^2(D)$ is a Hilbert subspace of $L^2(\partial
D,\sigma)$ with each element $f\in H^2(D)$ having a natural
holomorphic extension to $D$ (see \cite[Ch.~1.5]{krantz}). 

In the sequel, no 
generality is lost by taking $B$ in the statement of
Theorem~\ref{mainresult} to be the unit ball $B_1$, and
  the smaller concentric ball
to be $B_r$, the ball of radius $r$ centred at $0$.  

If $L:X_1\to X_2$ is a continuous operator
between Banach spaces then for $k\ge1$, its \emph{$k$-th approximation
  number} $a_k(L)$ is defined as
$$
a_k(L)= \inf \all{\norm{L-K}{}}{K:X_1\to X_2\text{ linear and
    continuous with }{\rm
    rank}(K)<k}\,.
$$
\end{notation}

The proof of Theorem \ref{mainresult} hinges on the following
two lemmas.

\begin{lemma}
\label{weylforL}
If $J:H^2(B_1)\hookrightarrow H^\infty(B_r)$ denotes the canonical embedding,
then $J$ and $\l$ are compact and for all $n\ge1$ 
\begin{equation}
\label{weylhilbertforL}
|\lambda_n(\l)| \le W  \prod_{k=1}^n a_k(J)^{1/n}\,.
\end{equation}
\end{lemma}
\begin{proof}
If $f\in H^2(B_1)$ and $z\in B_r$ then $|f(z)|\leq
(2/(1-r))^{d/2}$ by \cite[Thm.~7.2.5]{rudin}, so
$\all{f}{\norm{f}{H^2(B_1)}\leq 1}$ is a normal family in $H^\infty(B_r)$,
hence relatively compact in $H^\infty(B_r)$ by Montel's Theorem 
(see \cite[Ch.~1, Prop.~6]{narasimhan}), thus $J$ is
compact. 

Next observe that if $f\in H^\infty(B_1)$ then $f\in H^2(B_1)$ by
\cite[Thm.~5.6.8]{rudin} and the canonical embedding
$\hat{J}:H^\infty(B_1)\hookrightarrow H^2(B_1)$ is continuous of norm $1$,
because $\sigma(\partial B_1)=1$.
We claim that $\hat{\l}f:= \sum_{i\in\I} w_i\cdot f\circ T_i$ 
defines a continuous operator $\hat{\l}:H^\infty(B_r)\to
 H^\infty(B_1)$. To see this, fix $f\in H^\infty(B_r)$ 
and note that $w_i\cdot f\circ T_i\in H^\infty(B_1)$ 
with
$\norm{w_i\cdot f\circ T_i}{H^\infty(B_1)}\leq \norm{w_i}{H^\infty(B_1)}\norm{f}{H^\infty(B_r)}$ 
for every $i\in\I$. But since     
$\|\hat{\l}f\|_{H^\infty(B_1)}\leq
\sum_{i\in\I}\norm{w_i}{H^\infty(B_1)}\norm{f}{H^\infty(B_r)}$ and
$\sum_{i\in\I}\norm{w_i}{H^\infty(B_1)}<\infty$ by hypothesis, we conclude that
$\hat{\l} f\in H^\infty(B_1)$ and 
that $\hat{\l}$ is continuous. 
Now $|f(T_i(z))|\leq \norm{f}{H^\infty(B_r)}$ for every $z\in B_1$,
 $i\in\I$, so 
$\|\hat{\l}f\|_{H^\infty(B_1)}=\sup_{z\in B_1}|(\hat{\l}f)(z)|\leq 
\sup_{z\in B_1}
  \sum_{i\in\I}\abs{w_i(z)}\,\abs{f(T_i(z))}\leq
 W\norm{f}{H^\infty(B_r)}$,
and hence $\|\hat{\l}\|\le W$.
Now clearly 
$\l=\hat{J}\hat{\l}J$, so $\l$ is compact, and 
\begin{equation}
\label{embeddinglift}
a_k(\l)\le \|\hat{J}\hat{\l}\|a_k(J) \le Wa_k(J)
\quad\text{for all }k\ge1\,,
\end{equation}
since in general $a_k(L_1L_2)\le \norm{L_1}{} a_k(L_2)$ whenever $L_1$
and $L_2$ are bounded operators between Banach spaces (see
\cite[2.2]{pietsch}).
Moreover, since $\l$ is a
compact operator on Hilbert space,
Weyl's inequality (see \cite[3.5.1]{pietsch}, \cite{weyl})
asserts that
$ \prod_{k=1}^n \left| \lambda_k(\l)\right| \leq \prod_{k=1}^n
  a_k(\l)$
for all $n\ge1$.
Together with (\ref{embeddinglift}) this yields
(\ref{weylhilbertforL}),
because $|\lambda_n(\l)|\le \prod_{k=1}^n |\lambda_k(\l)|^{1/n}$.
\end{proof}

\begin{lemma}\label{hilbertembedding}
If $h_d(k):=\binom{k+d}{d}$
then for all $n\ge1$,
\begin{equation}
\label{an_seriesmodified}
a_n(J)^2\leq \sum_{l=k}^\infty h_{d-1}(l) r^{2l}\quad
\text{where $k\ge0$ is such that  }\  
h_d(k-1)<n\leq h_d(k)\,.
\end{equation}
\end{lemma}
\begin{proof}
$H^2(B_1)$ has reproducing
  kernel
$
  K(z,\zeta)=(1-(z,\zeta)_{\C^d})^{-d}
  $
(see \cite[Thm.~1.5.5]{krantz}\footnote{Note that the extra factor
  $(d-1)!/(2\pi^d)$ 
appearing in 
\cite[Thm.~1.5.5]{krantz}
is due to a different normalisation of the measure $\sigma$ on
$\partial B_1$.}),
where $(\cdot,\cdot)_{\C^d}$
denotes the Euclidean 
inner product, and
$K(z,\zeta)=\sum_{n=1}^\infty p_n(z)\overline{p_n(\zeta)}$
whenever
$\{p_n\}_{n=1}^\infty$ 
is an orthonormal basis for $H^2(B_1)$, 
the series converging pointwise
for every $(z,\zeta)\in B_1 \times B_1$ (see
\cite[p.~19]{halmos}).  

Define 
$J_n:H^2(B_1)\to H^\infty(B_r)$ by
$J_nf=\sum_{k=1}^{n-1} (f,p_k)\, p_k$.
If $z\in B_r$ then
\begin{multline*}
  |Jf(z)-J_nf(z)|^2 = |f(z) - J_nf(z)|^2 
     = \left| \sum_{k=n}^\infty
       (f,p_k)\, p_k(z)\right|^2 \\
    \le \sum_{k=n}^\infty
       |(f,p_k)|^2  \sum_{k=n}^\infty
    |p_k(z)|^2 
    \le \norm{f}{H^2(B_1)}^2 (K(z,z)-\sum_{k=1}^{n-1}|p_k(z)|^2)\,,
\end{multline*}
so 
\begin{equation}
\label{approxnumberformulainlemma}
a_n(J)^2 \le \sup_{z\in B_r}
\left( K(z,z)-\sum_{k=1}^{n-1}|p_k(z)|^2
\right)\,.
\end{equation}

If $n=1$ then $k=0$, 
in which case (\ref{an_seriesmodified}) follows from
(\ref{approxnumberformulainlemma})
since
$\sum_{l=0}^\infty h_{d-1}(l) r^{2l}=(1-r^2)^{-d}$.
Now define the
orthonormal basis
$\all{ p_{\underline{n}}}{\underline{n}\in\N_0^d}$ 
by (cf.~\cite[Prop.~1.4.8, 1.4.9]{rudin})
\begin{equation*}
p_{\underline{n}}(z)=K_{\underline{n}}z^{\underline{n}}
\quad (\underline{n}\in \N_0^d)\,,
\end{equation*}
where
$ 
K_{\underline{n}}
= \sqrt{\frac{(|\underline{n}|+d-1)!}{(d-1)!\,\underline{n}!}}$,  
$\un=(n_1,\ldots,n_d)$, $z^{\un}=z_1^{n_1}\cdots z_d^{n_d}$,
 $\un!=n_1!\cdots n_d!$, $|\un|=n_1+\cdots +n_d$.

If $n\ge2$ then 
there are $\binom{k+d-1}{d}$ multinomials
of degree less than or equal to $k-1$, so
\begin{equation*}
  a_n(J)^2 \leq \sup_{z\in B_r}\left ( K(z,z)-\sum_{\abs{\un}\leq
  k-1}\abs{p_{\un}(z)}^2 \right )
 =\sup_{z\in
  B_r}\sum_{l=k}^\infty\sum_{\abs{\un}=l}\abs{p_{\un}(z)}^2 
\leq \sum_{l=k}^\infty \frac{(l+d-1)!}{(d-1)!\,l!}r^{2l}
\end{equation*}
for all $n>\binom{k+d-1}{d}$,
because 
$\sum_{\abs{\un}=l}\frac{1}{\un !}\abs{z^{\un}}^2\leq
\frac{1}{l!}r^{2l}$
for $z\in B_r$ by the multinomial theorem. 
\end{proof}

\begin{proof}[Proof of Theorem~\ref{mainresult}]
By Lemma \ref{weylforL} it suffices to bound the geometric means
$(\prod_{k=1}^n a_k)^{1/n}$,
where $a_k:=a_k(J)$.
From Lemma \ref{hilbertembedding}
it follows that
\begin{equation}
\label{basiccorobound}
a_n^2 \leq \at_n\frac{r^{2\bt_n}}{(1-r^2)^d}\quad\text{for all }n\ge1\,,
\end{equation}
where 
\[ \begin{array}{l} \at_n:=h_{d-1}(k) \\ \bt_n:=k \end{array} \quad \text{
  for $h_{d}(k-1)<n\leq h_d(k)$}\,,\]
because
\begin{equation*}
\sum_{l=k}^\infty
h_{d-1}(l)r^{2l}  = h_{d-1}(k)r^{2k}\sum_{l=0}^\infty
\frac{h_{d-1}(l+k)}{h_{d-1}(k)}r^{2l} 
 \leq h_{d-1}(k)r^{2k} \sum_{l=0}^\infty h_{d-1}(l)r^{2l}
 =h_{d-1}(k)\frac{r^{2k}}{(1-r^2)^{d}}\,.
\end{equation*}

Combining (\ref{basiccorobound}) with Lemma \ref{weylforL}
gives, for all $n\ge1$,
\begin{equation}
\label{evalueW}
 |\lambda_n(\l)|\leq W \alpha_n\frac{r^{\beta_n}}{(1-r^2)^{d/2}}\,,
\end{equation}
where 
\[ \alpha_n :=\prod_{l=1}^n \at_l^{1/(2n)}\,,
\quad \beta_n:=\frac{1}{n}\sum_{l=1}^n\bt_l\,.
\]

To obtain (\ref{generalevalue}) and (\ref{d=1evalue})
from (\ref{evalueW})
we require an upper bound on $\alpha_n$ and a lower bound
on $\beta_n$. We start with the bounds for $\alpha_n$. 
Observe that 
\begin{equation}
\label{alphabounds}
\tilde \alpha_1=h_{d-1}(0)=1\,, \text{ and }
\tilde\alpha_l \le d(l-1)^{1-1/d}\ \text{for }l\ge2\,.
\end{equation}
To see this note that 
\begin{equation*} 
\frac{h_{d-1}(k)}{h_d(k-1)^{1-1/d}} =
\frac{(d!)^{1-1/d}}{(d-1)!}\left (
  \frac{\prod_{l=1}^{d-1}(k+l)^d}{\prod_{l=0}^{d-1}(k+l)^{d-1}}\right
)^{1/d}
=\frac{(d!)^{1-1/d}}{(d-1)!}\prod_{l=1}^{d-1}\left
  (1+\frac{l}{k}\right )^{1/d}
\end{equation*}
is decreasing in $k$, so if $h_d(k-1)<n\leq h_d(k)$ then
$
\frac{\at_l}{(l-1)^{1-1/d}}\leq 
\frac{h_{d-1}(k)}{h_d(k-1)^{1-1/d}}
\leq
\frac{h_{d-1}(1)}{h_d(0)^{1-1/d}}=d\,.
$

The estimate (\ref{alphabounds}) now yields the upper bound 
\begin{equation}
\label{alphanbound}
\alpha_n = \prod_{i=1}^n\at_i^{1/(2n)}
            \leq \sqrt{d} ((n-1)!)^{(d-1)/(2dn)} 
            \leq \sqrt{d} \left ( 2\left ( \frac{n}{e}\right )^n 
    \right )^{(d-1)/(2dn)}
\leq \sqrt{d}n^{(d-1)/(2d)} \,,
\end{equation}
where, for $n>1$, we have used the estimate
$(n-1)!\leq 2\left ( \frac{n}{e}\right )^n$
(i.e.~$ \log (n-1)!\leq \int_{x=2}^n \log x\, dx 
%= n\ln n-n -2\ln 2+2 
\leq n\log n-n+\log 2$).

We now turn to the bounds for $\beta_n$. 
If $h_d(k-1)<l\le h_d(k)$, so that $\tilde \beta_l=k$, then
$l\le h_d(k)\le (d!)^{-1}(k+d)^d$, which implies
$
\tilde \beta_l = k \ge (d!)^{1/d} l^{1/d} -d\,.
$
Therefore
\begin{equation}
\label{betanbound}
\beta_n = \frac{1}{n} \sum_{l=1}^n \tilde\beta_l
\ge -d +  (d!)^{1/d} \frac{1}{n} \sum_{l=1}^n l^{1/d} 
>  -d +  (d!)^{1/d} \frac{d}{d+1}n^{1/d}\,,
\end{equation}
where we have used 
$\sum_{l=1}^n l^{1/d} >  \int_{x=0}^n x^{1/d}\, dx =\frac{d}{d+1}
n^{1+1/d}$.
 
Assertion (\ref{generalevalue}) now follows from
(\ref{evalueW}), (\ref{alphanbound}), and (\ref{betanbound}). 
Finally, if $d=1$ then $\beta_n= \frac{1}{n} \sum_{l=1}^n \tilde\beta_l
=\frac{1}{n} \sum_{l=1}^n (l-1) = (n-1)/2$,
and (\ref{alphanbound}) becomes $\alpha_n\le 1$, so substituting into
(\ref{evalueW}) yields (\ref{d=1evalue}).
\end{proof}

\end{document}